\documentclass[A4paper,14pt]{article}

\usepackage{inputenc} 
\usepackage[T1]{fontenc}
\usepackage{amssymb,amsfonts,amsmath}
\newtheorem{definition}{Definition}

\newtheorem{proposition}[definition]{Proposition}
\newtheorem{corollary}[definition]{Corollary}
\newtheorem{fact}[definition]{Fact}
\newtheorem{lemma}[definition]{Lemma}
\newtheorem{theorem}[definition]{Theorem}

\newtheorem{problem}[definition]{Problem}

\newtheorem{claim}[definition]{Claim}

\newtheorem{conjecture}{Conjecture}

\newcommand{\bds}{\begin{description}}
\newcommand{\eds}{\end{description}}
\newcommand{\be}{\begin{enumerate}}
\newcommand{\ee}{\end{enumerate}}
\newcommand{\bq}{\begin{quote}}
\newcommand{\eq}{\end{quote}}

\newcommand{\qed}{\begin{flushright}{$\Box$}\end{flushright}}

\newcommand{\proof}{\noindent {\bf Proof: }}

\newcommand{\la}{\langle}
\newcommand{\ra}{\rangle}

\renewcommand{\qed}{\hfill $\square$\medskip}
\newcommand{\qedd}{\hfill $\blacksquare$\medskip}
\newcommand{\nothing}[1]{}

\def\cal{\mathcal}
\newcommand{\MM}{{\sf MM}}

\newcommand{\PFA}{{\sf PFA}}
\newcommand{\BPFA}{{\sf BPFA}}

\newcommand{\CH}{{\sf CH}}
\newcommand{\ZFC}{{\sf ZFC}}

\newcommand{\SCH}{\sf SCH}

\newcommand{\CP}{\sf CP}

\newcommand{\MRP}{{\sf MRP}}

\newcommand{\PID}{{\sf PID}}

\newcommand{\w}{\omega}
\renewcommand{\a}{\alpha}

\newcommand{\cc}{{\mathfrak c}}
\newcommand{\bb}{{\mathfrak b}}

\newcommand{\cf}{{\mbox{\tt cof}}}
\newcommand{\otp}{{\mbox{\tt otp}}}

\newcommand{\D}{{\mathcal D}}
\newcommand{\E}{{\mathcal E}}
\newcommand{\I}{{\mathcal I}}

\newcommand{\C}{\mathcal C}

\newcommand{\W}{{\mathcal W}}
\newcommand{\U}{{\mathcal U}}
\begin{document}
\title{\huge A family of covering properties for forcing axioms and strongly compact cardinals}
\author{Matteo Viale\\
KGRC, University of Vienna\\
Equipe de Logique Math\'ematique, Universit\'e Paris 7\\
{\texttt matteo@logic.univie.ac.at}}

\date{}
\maketitle
\begin{abstract}
\noindent This\footnote{This research has been partially supported
by the Austrian Science Fund FWF project P19375-N18. I thank Boban
Velickovic for many useful suggestions and comments on the redaction
of this paper.} paper presents the main results in my Ph.D. thesis
and develops from \cite{viaSCH} where it was shown that $\SCH$
follows from $\PFA$. In what follows several other and simpler
proofs of $\SCH$ are presented introducing a family of covering
properties $\CP(\kappa,\lambda)$ which implies both $\SCH$ and the
failure of square. I will also apply these covering properties to
investigate models of strongly compact cardinals or of strong
forcing axioms like $\MM$ or $\PFA$.
\end{abstract}

\noindent In this paper I introduce a family of covering properties
$\CP(\kappa,\lambda)$ indexed by pairs of regular cardinals
$\lambda<\kappa$. In the first part I show that these covering
properties capture the combinatorial content of many of the known
proofs of the failure of square-like principles from forcing axioms
or large cardinals. In the second part I show that a large class of
these covering properties follows either from the existence of a
strongly compact cardinal or from at least two combinatorial
principles which hold under $\PFA$ and are mutually independent: the
mapping reflection principle $\MRP$ introduced by Moore in
\cite{mooMRP} and the $P$-ideal dichotomy $\PID$ introduced in its
full generality by Todor\v{c}evi\'c in \cite{todPID} developing on
preceding works by him and Abraham \cite{abrtodPID} on $P$-ideals of
countable subsets of $\w_1$. This allows for a unified and simple
proof of the failure of square and of the singular cardinal
hypothesis $\SCH$ assuming $\PFA$. Finally in the last part of the
paper I will investigate the "saturation properties"\footnote{In a
sense that will be made explicit in the last section of the paper.}
of models of $\CP(\kappa,\lambda)$ for various $\kappa$ and
$\lambda$.

\noindent The paper is organized as follows: in sections
\ref{secPID} and \ref{secMRP} we introduce the combinatorial
principles $\PID$ and $\MRP$. In sections \ref{secCP} and
\ref{secCPBF} we introduce the covering properties
$\CP(\kappa,\lambda)$ and we outline some of their consequences
among which $\SCH$ and the failure of square. In sections
\ref{secCPLC}, \ref{secCPPID} and \ref{secCPMRP} we prove various
instances of $\CP(\kappa,\lambda)$ assuming respectively the
existence of a strongly compact cardinal, $\PID$ or
$\MRP$\footnote{The shortest path to obtain a self-contained proof
of the $\SCH$ starting from $\PID$ is to read sections \ref{secPID},
\ref{secCP} and \ref{secCPPID}, from $\MRP$ is to read sections
\ref{secMRP}, \ref{secCP} and \ref{secCPMRP}.}. In the  last section
we study the rigidity of models of $\CP(\kappa,\lambda)$ and we
prove several results that strongly suggest that any two models
$W\subseteq V$ of $\MM$ with the same cardinals have the same
$\w_1$-sequences of ordinals. Another result that we will present is
that strongly compact cardinals $\lambda$ are destroyed by any
forcing which preserves cardinals and changes the cofinality of some
regular $\kappa>\lambda$ to some $\theta<\lambda$.

\noindent The paper aims to be accessible and self-contained for any
reader with a strong background in combinatorial set theory.
Moreover while forcing axioms are a source of inspiration for the
results that we present, all the technical arguments in this paper
(except some of those in the last section) can be followed by a
reader with no familiarity with the forcing techniques. When not
otherwise explicitly stated \cite{jecST} is the standard source for
notation and definitions. For a regular cardinal $\theta$, we use
$H(\theta)$ to denote the structure $\la H(\theta),\in,<\ra$ whose
domain is the collection of sets whose transitive closure is of size
less than $\theta$ and where $<$ is a predicate for a fixed well
ordering of $H(\theta)$. For cardinals $\kappa\geq\lambda$ we let
$[\kappa]^\lambda$ be the family of subsets of $\kappa$ of size
$\lambda$. In a similar fashion we define $[\kappa]^{<\lambda}$,
$[\kappa]^{\leq\lambda}$, $[X]^\lambda$, where $X$ is an arbitrary
set. If $X$ is an uncountable set and $\theta$ a regular cardinal,
$\E\subseteq[X]^\theta$ is unbounded if for every $Z\in [X]^\theta$,
there is $Y\in\E$ containing $Z$. $\E$ is bounded otherwise. $\E$ is
closed in $[X]^\w$ if whenever $X=\bigcup_n X_n$ and $X_n\subseteq
X_{n+1}$ are in $\E$ for all $n$, then also $X\in\E$. It is a well
known fact that $\C\subseteq[X]^\omega$ is closed and unbounded
(club) iff there is $f:[X]^{<\omega}\rightarrow X$ such that $\C$
contains the set of all $Y\in [X]^{\omega}$ such that
$f[Y]^{<\omega}\subseteq Y$. $S\subseteq[X]^\omega$ is stationary if
it intersects all club subsets of $[X]^{\omega}$. The $f$-closure of
$X$ is the smallest $Y$ containing $X$ such that
$f[[Y]^{<\omega}]\subseteq Y$. Given $f$ as above, $\E_f$ is the
club of $Z\in[X]^\omega$ such that $Z$ is $f$-closed. If $X$ is a
set of ordinals, then $\overline{X}$ denotes the topological closure
of $X$ in the order topology. For regular cardinals
$\lambda<\kappa$, $S_\kappa^{\leq\lambda}$ denotes the subset of
$\kappa$ of points of cofinality $\leq\lambda$. In a similar fashion
we define $S_\kappa^\lambda$ and $S_\kappa^{<\lambda}$. We say that
a family $\mathcal{D}$ is covered by a family $\mathcal{E}$ if for
every $X\in\mathcal{D}$ there is a $Y\in\mathcal{E}$ such that
$X\subseteq Y$. We also recall the following definitions central to
the arguments that follows:

\smallskip

\noindent The singular cardinal hypothesis $\SCH$ asserts that
$\kappa^{\cf(\kappa)}=\kappa^++2^{\cf(\kappa)}$ for all infinite
cardinals $\kappa$.

\smallskip

\noindent It is a celebrated result of Silver \cite{sil75} that if
$\SCH$ fails, then it first fails at a singular cardinal of
countable cofinality. It is also known that the failure of $\SCH$ is
a strong hypothesis as it entails the existence of models of $\ZFC$
with measurable cardinals \cite{gitFSCH}.

\smallskip

\noindent Let $\kappa$ be an infinite regular cardinal. The square
principle $\Box(\kappa)$ asserts the existence of a sequence
$(C_\a:\a<\kappa)$ with the following properties:
\begin{itemize}
\item[\it (i)] for every limit $\a$, $C_\a$ is a closed unbounded
subset of $\a$,
\item[\it (ii)] if $\a$ is a limit point of $C_\beta$,
$C_\a=C_\beta\cap\a$,
\item[\it (iii)] there is no club $C$ in $\kappa$ such that for all
$\a$ there is $\beta\geq\a$ such that $C\cap\a= C_\beta$,
\item[\it (iv)] $C_{\beta+1}=\{\beta\}$.
\end{itemize}

\noindent It is well known that the failure of $\Box(\kappa)$ for
all cardinals $\kappa>\aleph_1$ is a very strong large cardinal
hypothesis. For example Schimmerling has shown that it entails the
existence of models of set theory with Woodin cardinals
\cite{schBOX}.

\smallskip

\noindent Recall that $\lambda$ is a strongly compact cardinal if
for every $\kappa\geq\lambda$ there is a $\lambda$-complete
ultrafilter on $[\kappa]^{<\lambda}$.

\section{The $P$-ideal dichotomy} \label{secPID}

\noindent Let $Z$ be an uncountable set. $\I\subseteq[Z]^{\leq\w}$
is a $P$-ideal if it is an ideal and for every countable family
$\{X_n\}_n\subseteq\I$ there is an $X\in\I$ such that for all $n$,
$X_n\subseteq^*X$ (where $\subseteq^*$ is inclusion modulo finite).

\begin{definition}(Todor\v{c}evi\'c, \cite{todPID})

\noindent The $P$-ideal dichotomy ($\PID$) asserts that for every
$P$-ideal $\I$ on $[Z]^{\leq\w}$ for some fixed uncountable $Z$, one
of the following holds:

\begin{itemize}

\item[\it (i)] There is $Y$ uncountable subset of $Z$ such that
$[Y]^{\leq\w}\subseteq\I$.

\item[\it (ii)] $Z=\bigcup_n A_n$ with the property that
$A_n$ is orthogonal to $\I$ (i.e. $X\cap Y$ is finite for all
$X\in[A_n]^\w$ and $Y\in\I$) for all $n$.

\end{itemize}

\end{definition}

\noindent $\PID$ is a principle which follows from $\PFA$ and which
is strong enough to rule out many of the standard consequences of
$V=L$. For example Abraham and Todor\v{c}evi\'c \cite{abrtodPID}
have shown that under $\PID$ there are no Souslin trees while
Todor\v{c}evi\'c has shown that $\PID$ implies the failure of
$\square(\kappa)$ on all regular $\kappa>\aleph_1$ \cite{todPID}.
Due to this latter fact the consistency strength of this principle
is considerable. Another interesting result by Todor\v{c}evi\'c is
that $\PID$ implies that $\bb\leq\aleph_2$\footnote{In \cite{todPID}
it is shown that any gap in $P(\w)/FIN$ is either an Haussdorff gap
or a $(\kappa,\w)$ gap with $\kappa$ regular and uncountable. By
another result of Todor\v{c}evi\'c (see \cite{jecST} pp. 578 for a
proof) if $\bb>\aleph_2$ there is an $(\w_2,\lambda)$ gap in
$P(\w)/FIN$ for some regular uncountable $\lambda$. Thus $\PID$ is
not compatible with $\bb>\aleph_2$.}. Nonetheless in
\cite{abrtodPID} and \cite{todPID} it is shown that this principle
is consistent with $\CH$. Other interesting applications of $\PID$
can be found in \cite{velPID}, \cite{baljecpazPIDCCC}, \cite{todPID}
and \cite{abrtodPID}.

\section{The mapping reflection principle} \label{secMRP}

Almost all known applications of $\MM$ which do not follow from
$\PFA$ are a consequence of some form of reflection for stationary
sets. These types of reflection principles are a fundamental source
in order to obtain proofs of all cardinal arithmetic result that
follows from $\MM$. In particular $\SCH$ and the fact that
$\cc\leq\w_2$ are a consequence of many of the known reflection
principles which hold under $\MM$. However up to a very recent time
there was no such kind of principle which could be derived from
$\PFA$ alone. This has been the main difficulty in the search for a
proof that $\PFA$ implies $\cc=\w_2$, a result which has been
obtained by Todor\v{c}evi\'c and Veli\v{c}kovi\'c appealing to
combinatorial arguments which are not dissimilar from the $P$-ideal
dichotomy \cite{velSCH}. Later on this has also been the crucial
obstacle in the search for a proof of $\SCH$ from $\PFA$.

\noindent In $2003$ Moore \cite{mooMRP} found an interesting form of
reflection which can be derived from $\PFA$, the mapping reflection
principle $\MRP$. He has then used this principle to show that
$\BPFA$ implies that $\cc=\aleph_2$ and also that this principle is
strong enough to entail the non-existence of square sequences. He
has also shown in \cite{mooSCH} that $\MRP$ could be a useful tool
in the search of a proof of $\SCH$ from $\PFA$. I first obtained my
proof of this latter theorem elaborating from \cite{mooSCH}. Many
other interesting consequences of this reflection principle have
been found by Moore and others. A complete presentation of this
subject will be found in \cite{caivelMRP}.

\begin{definition}
Let $\theta$ be a regular cardinal, let $X$ be uncountable, and let
$M\prec H(\theta)$ be countable such that $[X]^{\omega}\in M$. A
subset $\Sigma$ of $[X]^{\omega}$ is $M$-stationary if for all
$\E\in M$ such that $\E\subseteq [X]^{\omega}$ is club, $\Sigma \cap
\E\cap M \neq \emptyset$.
\end{definition}

\noindent Recall that the Ellentuck topology on $[X]^{\omega}$ is
obtained by declaring a set open if it is the union of sets of the
form
$$
[x,N]=\{ Y\in [X]^{\omega} \colon x\subseteq Y\subseteq N \}
$$
\noindent where $N\in[X]^\w$ and $x\subseteq N$ is finite.

\begin{definition}
$\Sigma$ is an open stationary set mapping if there is an
uncountable set $X$ and a regular cardinal $\theta$ such that
$[X]^{\omega}\in H(\theta)$, the domain of $\Sigma$ is a club in
$[H(\theta)]^{\omega}$ of countable elementary submodels $M$ such
that $X\in M$ and for all $M$, $\Sigma(M)\subseteq [X]^{\omega}$ is
open in the Ellentuck topology on $[X]^\omega$ and $M$-stationary.
\end{definition}

\noindent The mapping reflection principle ($\MRP$) asserts that:
\begin{quote}
If $\Sigma$ is an open stationary set mapping, there is a continuous
$\in$-chain $\vec N= (N_\xi\colon\xi<\w_1)$ of elements in the
domain of $\Sigma$ such that for all limit ordinals $0<\xi<\w_1$
there is  $\nu<\xi$ such that $N_\eta\cap X\in\Sigma(N_\xi)$ for all
$\eta$ such that $\nu <\eta <\xi$.
\end{quote}
\noindent If $(N_{\xi}\colon \xi <\omega_1)$ satisfies the
conclusion of $\MRP$ for $\Sigma$ then it is said to be a reflecting
sequence for $\Sigma$.

\smallskip

\noindent We list below some of the interesting consequences of
$\MRP$.

\begin{itemize}

\item (Moore \cite{mooMRP}) $\PFA$ implies $\MRP$.

\item (Moore \cite{mooMRP}) $\MRP$ implies that $\cc=\aleph_2$.
As a simple outcome of his proof of the above theorem Moore obtains
also that $\BPFA$ implies $\cc=\aleph_2$.

\item (Moore \cite{mooMRP}) \label{thMRPBOX} Assume $\MRP$. Then
$\Box(\kappa)$ fails for all regular $\kappa\geq\aleph_2$.

\end{itemize}

\noindent A folklore problem in combinatorial set theory for the
last twenty years has been the consistency of the existence of a
five element basis for the uncountable linear orders, i.e. the
statement that there are five uncountable linear orders such that at
least one of them embeds in any other uncountable linear order.

\begin{itemize}

\item \label{thMRPSC} (Moore \cite{mooMRPSC}) Assume $\BPFA$ and
$\MRP$. Then there is a five element basis for the uncountable
linear orders.

\end{itemize}

\noindent A considerable reduction of the large cardinal hypothesis
needed for the consistency of the above conjecture has been obtained
in \cite{konlarmoovelMRP}. A byproduct of their results yields to
the following:

\begin{itemize}
\item (K\"onig, Larson, Moore, Veli\v{c}kovi\'c
\cite{konlarmoovelMRP}) $\MRP$ implies that there are no Kurepa
trees.
\end{itemize}

\noindent Other interesting consequences of $\MRP$ can be found in
\cite{caivelBPFA}.

\noindent We remark that $\MRP$ and $\PID$ are mutually independent
principles since $\PID$ is compatible with $\CH$ while, by a result
of Myiamoto \cite{miyMRPSH}, $\MRP$ is compatible with the existence
of Souslin trees.

\section{A family of covering properties}\label{secCP}

\noindent In this section we introduce the main original concept of
this paper. A careful analysis of Todor\v{c}evi\'c's proof that
$\PID$ implies that $\Box(\kappa)$ fails for all regular
$\kappa\geq\aleph_1$ leads to the isolation of a family of covering
properties $\CP(\kappa,\lambda)$ which on one side are strong enough
to entail both $\SCH$ and the failure of square, and on the other
side are weak enough to be a consequence of the existence of a
strongly compact cardinal, of $\PID$ and of $\MRP$.

\begin{definition} \label{defCP}
For regular cardinals $\lambda<\kappa$, $\D$
$=(K(\a,\beta):\a<\lambda,\, \beta\in \kappa)$ is a
$\lambda$-covering matrix for $\kappa$ if: \bds
\item[\it (i)] $\beta\subseteq\bigcup_{\a<\lambda} K(\a,\beta)$ for all $\beta$,
\item[\it (ii)] $|K(\a,\beta)|<\kappa$ for all $\beta$ and $\a$,
\item[\it (iii)] $K(\a,\beta)\subseteq K(\eta,\beta)$ for all $\beta<\kappa$ and for
all $\a<\eta<\lambda$,
\item[\it (iv)]  for all $\gamma<\beta<\kappa$ and for all $\a<\lambda$, there is
$\eta<\lambda$ such that $K(\a,\gamma)\subseteq K(\eta,\beta)$.
\item[\it (v)] for all $X\in[\kappa]^{\leq\lambda}$, there is $\gamma_X<\kappa$
such that for all $\beta<\kappa$ and $\eta<\lambda$, there is $\a$
such that $K(\eta,\beta)\cap X\subseteq K(\a,\gamma_X)$ \eds

\noindent $\beta_\D\leq\kappa$ is the least $\beta$ such that for
all $\a$ and $\gamma$, $\otp(K(\a,\gamma))<\beta$
\end{definition}

\noindent We will mainly be interested in $\w$-covering matrices,
which we will just call covering matrices. As we will see below
square like principles are useful to construct several kinds of
covering matrices. One successful strategy to negate square
principles from large cardinals and forcing axioms is to use
appropriate ultrafilters or specific forcing arguments to
"diagonalize" through the covering matrix defined appealing to these
square-like principles, as for example in the proofs of the failure
of square from a strongly compact by Solovay \cite{solSCH} or from
$\PID$ by Todor\v{c}evi\`c \cite{todPID}\footnote{See also the
several arguments of this sort appearing in the sequel of this
paper.}. The covering matrices induced by square-like principles
that we will consider satisfy a stronger coherence property than the
"local" property {\it (v)}. This condition is replaced by the
"global" property\footnote{For example the matrix produced by a
square sequence using the $\rho_2$-function (see sections $6$ and
$8$ of \cite{todHST} and fact \ref{faCPBOX} below) satisfies {\it
(i)},$\ldots$ {\it (iv)} and {\it (v')}, and this matrix can be used
to show that $\CP$ implies the failure of square. Another
interesting example of a covering matrix which is not defined
appealing to lemma \ref{lemCP} below and which satisfies {\it (v)}
is the matrix used in the proof of theorem \ref{thPFCP} in the last
section.}:

\begin{quote}
{\it (v') For all $\gamma<\beta<\kappa$ and $\eta<\lambda$, there is
$\a$ such that $K(\eta,\beta)\cap \gamma\subseteq K(\a,\gamma)$.}
\end{quote}

\noindent The key point to introduce condition {\it (v)} above in
this weak form is that $\lambda$-covering matrices on $\kappa$ can
be defined in an elementary way in $\ZFC$ and the diagonalization
argument which in the square-like cases leads to a contradiction, in
the general case leads to a simple combinatorial argument to compute
$\kappa^\lambda$. It is possible to prove the following:

\begin{lemma}
Assume that there is a stationary set of points of uncountable
cofinality in the approachability ideal\footnote{In \cite{sheCA} it
is possible to find a definition of the approachability ideal
$\I[\kappa]$. We avoid it in this paper since it is not relevant for
the arguments we are presenting.} $\I[\kappa]$. Then there is an
$\w$-covering matrix on $\kappa$. Moreover if $\lambda$ is singular
of countable cofinality, then there is $\D$ covering matrix on
$\lambda^+$ with $\beta_\D=\lambda$.
\end{lemma}

\noindent Since the main original application of the existence of
$\w$-covering matrices $\D$ that we have found is the proof of
$\SCH$ from $\PFA$, we will just prove a weaker form of the lemma:

\begin{lemma} \label{lemCP}
Assume that $\kappa>\cc$ is a singular cardinal of countable
cofinality. Then there is an $\w$-covering matrix $\C$ for
$\kappa^+$ with $\beta_\C=\kappa$.
\end{lemma}

\proof The matrices we are going to define satisfy {\it
(i),$\cdots$,(iii)} and a stronger coherence property than what is
required by {\it (iv)} and {\it (v)} of the above definition. They
will satisfy the following properties {\it (iv*)} and {\it (v*)}
from which {\it (iv)} and {\it (v)} immediately follow:

\bds

\item[\it (iv*)] For all $\a<\beta$ there is $n$ such that $K(m,\a)\subseteq
K(m,\beta)$ for all $m\geq n$.

\item[\it (v*)] For all $X\in[\kappa^+]^\w$ there is $\gamma_X<\kappa$
such that for all $\beta\geq\gamma_X$ there is $n$ such that
$K(m,\beta)\cap X=K(m,\gamma_X)\cap X$ for all $m\geq n$.

\eds

\noindent Let $\phi_\eta:\kappa\rightarrow\eta$ be a surjection for
all $0<\eta<\kappa^+$. Fix also $\{\kappa_n:n<\w\}$ increasing
sequence of regular cardinals cofinal in $\kappa$ with
$\kappa_0\geq\aleph_1$. Set
$$
K(n,\beta)=\bigcup\{K(n,\gamma):\gamma\in\phi_\beta[\kappa_n]\}. $$
\noindent It is immediate to check that
$\D=(K(n,\beta):n\in\w,\beta<\kappa^+)$ satisfies {\it(i)},
{\it(ii)}, {\it(iii)} of definition \ref{defCP}, property {\it
(iv*)} and that $\beta_\D=\kappa$. To prove {\it (v*)}, let
$X\in[\kappa^+]^\w$ be arbitrary. Now since $\cc<\kappa^+$ and there
are at most $\cc$ many subsets of $X$, there is a stationary subset
$S$ of $\kappa^+$ and a fixed decomposition of $X$ as the increasing
union of sets $X_n$ such that $X\cap K(n,\a)=X_n$ for all $\a$ in
$S$ and for all $n$. Now properties {\it (i)}$\cdots${\it (iv)} of
the matrix guarantees that this property of $S$ is enough to get
{\it (v*)} for $X$ with $\gamma_X=\min(S)$. \qedd

\noindent A similar argument can be used to prove the following:

\begin{lemma} \label{lemCP1}
Assume that $\kappa$ and $\lambda$ are regular with
$\kappa>2^\lambda$. Then there is a $\lambda$-covering matrix $\C$
for $\kappa$.
\end{lemma}

\noindent We say that $\D$ is covered by $\E$ iff every $X\in\D$ is
contained in some $Y\in\E$.

\begin{definition} $\CP(\kappa,\lambda)$:
$\kappa$ has the $\lambda$-covering  property\footnote{Moore has
first noticed that a covering property similar to the $\w$-covering
property for a regular $\kappa>\cc$ followed from $\MRP$, reading a
draft of \cite{viaSCH}.} if for every $\D$, $\lambda$-covering
matrix for $\kappa$ there is an unbounded subset $A$ of $\kappa$
such that $[A]^{\lambda}$ is covered by $\D$. $\CP(\kappa)$
abbreviates $\CP(\kappa,\w)$ and $\CP$ is the statement that
$\CP(\kappa)$ holds for all regular $\kappa>\cc$.
\end{definition}

\begin{fact} \label{faCPSCH}
Assume $\CP(\kappa^+)$ for all singular $\kappa$ of countable
cofinality. Then $\lambda^{\aleph_0}=\lambda$, for every $\lambda
\geq 2^{\aleph_0}$ of uncountable cofinality.
\end{fact}

\noindent {\bf Proof:} By induction. The base case is trivial. If
$\lambda =\kappa ^+$ with $\cf(\kappa)> \omega$, then
$\lambda^{\aleph_0}=\lambda \cdot
\kappa^{\aleph_0}=\lambda\cdot\kappa=\lambda$, by the inductive
hypothesis on $\kappa$. If $\lambda$ is a limit cardinal and
$\cf(\lambda)>\omega$ then $\lambda^{\aleph_0}= \sup \{
\mu^{\aleph_0} : \mu <\lambda\}$, so the result also follows by the
inductive hypothesis. Thus, the only interesting case is when
$\lambda=\kappa^+$, with $\kappa$ singular of countable cofinality.
In this case we will show, using $\CP$, that
$(\kappa^+)^{\aleph_0}=\kappa^+$. To this aim let $\mathcal{D}$ be a
covering matrix for $\kappa^+$ with $\beta_{\D}=\kappa$. Notice that
by our inductive assumptions, since every $K(n,\beta)$ has order
type less than $\kappa$, $|[K(n,\beta)]^\w|$ has size less than
$\kappa$. So
$\bigcup\{[K(n,\beta)]^\w:n<\omega\,\&\,\beta\in\kappa^+\}$ has size
$\kappa^+$. Use $\CP$ to find $A\subseteq\kappa^+$ unbounded in
$\kappa^+$, such that $[A]^\omega$ is covered by $\mathcal{D}$. Then
$[A]^\omega\subseteq\bigcup\{[K(n,\beta)]^\omega:n<\omega\,\&\,\beta\in\kappa^+\}$,
from which the conclusion follows. \hfill$\blacksquare$\medskip

\smallskip

\noindent The following theorems motivate the introduction of these
covering properties:

\begin{theorem} \label{thSCSCH}
Assume $\lambda$ is strongly compact. Then $\CP(\kappa,\theta)$
holds for all regular $\theta<\lambda$ and all regular
$\kappa\geq\lambda$.
\end{theorem}

\begin{theorem} \label{thPIDSCH}
Assume $\PID$. Then $\CP$ holds.
\end{theorem}

\noindent On the other hand $\MRP$ allows us to infer a slightly
weaker conclusion than the one of the previous theorem.
\begin{theorem} \label{thMRPSCH}
Assume $\MRP$ and let $\D$ be a covering matrix for $\kappa$ such
that $K(n,\beta)$ is a closed set of ordinals for all $K(n,\beta)$.
Then there is $A$ unbounded in $\kappa$ such that $[A]^\w$ is
covered by $\D$.
\end{theorem}

\noindent In particular we obtain:

\begin{corollary}
$\PFA$ implies $\SCH$.
\end{corollary}

\proof $\PFA$ implies $\PID$ and $\PID$ implies $\CP$. In particular
$\PFA$ implies that $\kappa^\w=\kappa$ for all regular
$\kappa\geq\cc$. By Silver's theorem \cite{sil75} the least singular
$\kappa>2^{\cf\kappa}$ such that $\kappa^{\cf\kappa}>\kappa^+$ has
countable cofinality. Now assume $\PFA$ and let $\kappa$ have
countable cofinality. By fact \ref{faCPSCH},
$\kappa^{\cf(\kappa)}\leq(\kappa^+)^{\aleph_0}=\kappa^+$. Thus
assuming $\PFA$ there cannot be a singular cardinal of countable
cofinality which violates $\SCH$. Combining this fact with Silver's
result we get that $\SCH$ holds under $\PFA$. \qedd

\noindent Before proving all the above theorems we analyze in more
details the effects of $\CP$ and we give other interesting examples
of $\lambda$-covering matrices.

\section{Some simple consequences of $\CP(\kappa,\lambda)$} \label{secCPBF}

We remind the reader that if $\lambda^{\aleph_0}=\lambda$ for all
regular $\lambda>2^{\aleph_0}$, then by Silver's theorem and some
elementary cardinal arithmetic we obtain that
$\kappa^\lambda=\kappa+2^\lambda$ for all regular cardinals
$\lambda<\kappa$.

\begin{fact} \label{faCPTFAE}
Let $\theta<\kappa$ be regular cardinals such that
$\lambda^{\theta}<\kappa$ for all regular $\lambda<\kappa$.
$\D=\{K(\a,\beta):\a\in\theta,\beta<\kappa\}$ be a $\theta$-covering
matrix on $\kappa$ and $A$ be an unbounded subset of $\kappa$. The
following are equivalent:
\begin{itemize}
\item[\it(i)] $[A]^\theta$ is covered by $\D$.
\item[\it (ii)] $[A]^\lambda$ is covered by $\D$ for all $\lambda<\kappa$ such that
$\cf(\lambda)>\theta$.
\end{itemize}
\end{fact}

\proof \noindent {\it (ii)} implies {\it (i)} is evident. To prove
the other direction, assume {\it (i)} and let $Z\subseteq A$ have
size $\lambda$ with $\cf(\lambda)>\theta$. We need to find
$\a<\theta$ and $\beta<\kappa$ such that $Z\subseteq K(\a,\beta)$.
For $X\in [Z]^\theta\subseteq [A]^\theta$ let by {\it (i)}
$\a_X<\theta$, $\beta_X<\kappa$ be such that $X\subseteq
K(\a_X,\beta_X)$. By our assumptions, $\lambda^{\theta}<\kappa$. For
this reason $\beta=\sup_{X\in [Z]^\theta}\beta_X<\kappa$. Now by
property {\it (iii)} of $\mathcal{D}$, we have that for all $X\in
[Z]^\theta$, $X\subseteq K(\a_X,\beta)$ for some $\a_X$. Let $\C_\a$
be the set of $X$ such that $\a_X=\a$. Now notice that for at least
one $\a$, $\C_\a$ must be unbounded in $[Z]^\theta$, otherwise
$[Z]^\theta$ would be the union of $\theta$ bounded subsets which is
not possible since $|Z|$ has cofinality larger than $\theta$. Then
$Z\subseteq K(\a,\beta)$, since every $\alpha\in Z$ is in some
$X\in\C_\a$ because $\C_\a$ is unbounded. This completes the proof
of the fact. \qed




\noindent A $\lambda$-covering matrix $\D$ for $\kappa^+$ with
$\beta_\D=\kappa$ is an object simple to define when $\kappa$ has
cofinality $\lambda$. If $\kappa^\lambda=\kappa$, the existence of a
covering matrix $\D$ for $\kappa^+$ with $\beta_\D<\kappa^+$ is not
compatible with $\CP(\kappa^+,\lambda)$. This is a simple
consequence of the above facts:

\begin{fact} \label{faCMR}
Assume $\kappa^\lambda=\kappa$ and $\CP(\kappa^+,\lambda)$. Then
there is no $\lambda$-covering matrix $\D$ on $\kappa^+$ with
$\beta_\D<\kappa^+$.
\end{fact}

\proof Assume not and let $\D$ be a $\lambda$-covering matrix for
$\kappa^+$ with $\beta_\D<\kappa^+$. By $\CP(\kappa^+,\lambda)$
there should be an $A$ unbounded in $\kappa^+$ such that
$[A]^\lambda$ is covered by $\D$. Appealing to fact \ref{faCPTFAE},
we can conclude in any case that $[A]^\kappa$ is covered by $\D$.
Take $\beta$ large enough in order that
$\otp(A\cap\beta)>\beta_{\D}$. Since $A\cap\beta$ has size at most
$\kappa$ there are $n,\gamma$ such that $A\cap\beta\subseteq
K(n,\gamma)$. Thus $\beta_{\D}<\otp(A\cap\beta)\leq\otp
K(n,\beta)<\beta_\D$ a contradiction. \qedd

\noindent On the other hand assuming that $\kappa$ is regular to
what extent can we stretch $\lambda$ in order that
$\CP(\kappa^+,\lambda)$ is consistent? The fact below shows that
$\lambda$ cannot be $\kappa$.

\begin{lemma} \label{thCPmax}$\CP(\kappa^+,\kappa)$ fails for all regular
$\kappa\geq\w_1$.
\end{lemma}

\proof Fix a sequence $\C=\{C_\xi:\xi<\kappa^+\}$ such that for all
limit $\a$, $C_\a$ is a club subset of $\a$ of order type $\cf(\a)$.
If $\xi=\a+1$, $C_{\xi}=\{\a\}$. Define by induction on
$\a\leq\beta<\kappa^+$,
$$
\rho^*(\a,\beta):\,[\kappa^+]^2\rightarrow\kappa
$$
\noindent as follows:
\begin{itemize}
\item
$\rho^*(\a,\a)=0$,
\item
$\rho^*(\a,\beta)=\max\{\otp(C_\beta\cap\a),\rho^*(\a,\min(C_\beta\setminus\a)),\\
\sup\{\rho^*(\xi,\a):\xi\in C_\beta\cap\a\}\}$.
\end{itemize}

\noindent We will need the following properties of $\rho^*$ which
follows from the fact that $\kappa$ is regular\footnote{For a proof
see \cite{todHST} lemmas 19.1 and 19.2}.

\begin{lemma} \label{keylemrho1}
For all $\a\leq\beta\leq\gamma$:

\begin{description}

\item[(a)] $\rho^*(\a,\beta)\leq\max\{\rho^*(\a,\gamma),\rho^*(\beta,\gamma)\}$

\item[(b)] $\rho^*(\a,\gamma)\leq\max\{\rho^*(\a,\beta),\rho^*(\beta,\gamma)\}$

\end{description}
\qed
\end{lemma}

\begin{lemma}\label{keylemrho2}
For all $\a<\kappa^+$ and $\nu<\kappa$:
$$
|\{\xi<\a:\rho^*(\xi,\a)<\nu\}|\leq|\nu|+\aleph_0.
$$
\qed
\end{lemma}

\noindent Set $D(\a,\beta)=\{\xi<\beta:\rho^*(\xi,\beta)\leq\a\}$
for all $\a<\kappa$ and $\beta<\kappa^+$.

\begin{fact}\label{kayfacPFCP} The following holds:
\begin{description}

\item[\it (i)] $\otp(D(\a,\beta))<\kappa$ for all $\a<\kappa$ and $\beta<\kappa^+$,

\item[\it (ii)] for all $\gamma<\beta<\kappa^+$,there is $\a_0<\kappa$ such that
$D(\a,\gamma)=D(\a,\beta)\cap\gamma$ for all $\a\geq \a_0$.

\end{description}
\end{fact}

\proof {\it (i)} follows from the second lemma on $\rho^*$. To prove
{\it (ii)}, let $\a_0$ be such that $\rho^*(\gamma,\beta)\leq \a_0$,
$\a\geq \a_0$ and $\xi<\beta$ such that $\rho^*(\xi,\beta)\leq \a$ .
By {\it (a)} of lemma \ref{keylemrho1}:
$$
\rho^*(\xi,\gamma)\leq\max\{\rho^*(\xi,\beta),\rho^*(\gamma,\beta)\}
\leq\max\{\a,\a_0\}=\a.
$$
\noindent Conversely assume that $\rho^*(\xi,\gamma)\leq \a$, by
{\it (b)} of the same lemma:
$$
\rho^*(\xi,\beta)\leq\max\{\rho^*(\xi,\gamma),\rho^*(\gamma,\beta)\}
\leq\max\{\a,\a_0\}=\a.
$$

\noindent Thus $D(\a,\gamma)=D(\a,\beta)\cap\gamma$ for all $\a\geq
\a_0$. \qed

\noindent This means that
$\D=(D(\a,\beta):\a<\kappa,\beta<\kappa^+)$ is a $\kappa$-covering
matrix for $\kappa^+$ with $\beta_\D=\kappa$. Assuming
$\CP(\kappa^+,\kappa)$ there would be $A$ unbounded in $\kappa^+$
such that $[A]^\kappa$ is covered by $\D$. However there cannot be
an unbounded subset $A$ of $\kappa^+$ such that $[A]^\kappa$ is
covered by $\D$, since any element of the matrix has order type less
than $\kappa$. \qed

\noindent The function $\rho^*$ defined above will be useful also
for some other applications of $\CP(\kappa,\lambda)$ in the final
section (see theorem \ref{thPFCP}).

\noindent The following theorem follows closely Todor\v{c}evi\'c's
proof that $\PID$ entails the failure of square and shows that $\CP$
is a very large cardinal property.

\begin{theorem} \label{faCPBOX} Assume $\kappa>\aleph_1$ is regular.
Then $\CP(\kappa)$ implies that $\square(\kappa)$ fails.
\end{theorem}

\proof Todor\v{c}evi\'c has shown that assuming $\square(\kappa)$ it
is possible to define a step function (see sections $6$ and $8$ of
\cite{todHST}):
$$
\rho_2:\, [\kappa]^2\rightarrow\,\w
$$
\noindent with the following properties: \bds
\item[\it (i)] For every $A$ unbounded in $\kappa$, $\rho_2[[A]^2]$ is
unbounded in $\w$,
\item[\it (ii)] for every $\a<\beta$ there is $m$ such that
$|\rho_2(\xi, \a)-\rho_2(\xi,\beta)|\leq m$ for all $\xi<\a$. \eds

\noindent By {\it (ii)} it is immediate to check that
$\D=\{K(n,\a):n\in\w\,\&\,\a<\kappa\}$ is a covering matrix for
$\kappa$, where $K(n,\a)=\{\xi:\rho_2(\xi,\a)\leq n\}$. In fact it
can be shown something stronger i.e. that for every $\a<\beta$ and
$n$ there is $m$ such that $K(n,\a)\subseteq K(m,\beta)$ and
$K(n,\beta)\cap\a\subseteq K(m,\a)$.

\noindent Using this coherence property of $\D$ one gets that
whenever $A$ is an unbounded subset of $\kappa$ such that $[A]^\w$
is covered by $\D$, then for all $\beta<\kappa$,
$A\cap\beta\subseteq K(m_\beta,\beta)$ for some $m_\beta$. Thus one
can refine any such $A$ to an unbounded $B$ such that for a fixed
$m$, $B\cap\beta\subseteq K(m,\beta)$ for all $\beta\in B$. This
contradicts property {\it (i)} of $\rho_2$. Assuming $\CP(\kappa)$
we would get that an $A$ unbounded in $\kappa$ and such that
$[A]^\w$ is covered by $\D$ exists. However we just remarked that
this is impossible. \qed

\noindent The main difficulty towards a proof that $\PFA$ implies
$\SCH$ has been the fact that all standard principles of reflection
for stationary sets do not hold for $\PFA$. In particular Baudoin
\cite{bauPFASR} and Magidor (unpublished) have shown that $\PFA$ is
compatible with the existence on any regular $\kappa\geq\aleph_2$ of
a never reflecting stationary subset of $S_{\kappa}^\w$. However the
following form of reflection holds under $\CP(\kappa)$:

\begin{fact} \label{faREFPFA}
Assume $\CP(\kappa)$ for a regular $\kappa>\cc$ and let $\D$ be a
covering matrix for $\kappa$ with all $K(n,\beta)$ closed. Let
$\lambda<\kappa$ be a regular cardinal and let
$(S_\eta:\eta<\lambda)$ be an arbitrary family of stationary subsets
of $S_{\kappa}^{\leq\lambda}$. Then there exist $n$ and $\beta$ such
that $S_\eta\cap K(n,\beta)$ is non-empty for all $\eta<\lambda$.
\end{fact}

\proof By $\CP(\kappa)$ and facts \ref{faCPTFAE}, there is $X$
unbounded in $\kappa$ such that $[X]^\lambda$ is covered by $\D$.
Since $K(n,\beta)$ is closed for all $n$ and $\beta$, we have that
$[\overline{X}\cap S_{\kappa}^{\leq\lambda}]^{\lambda}$ is covered
by $\D$. To see this, let $Z$ be in this latter set and find
$Y\subseteq X$ of size $\lambda$ such that $Z\subseteq\overline{Y}$.
Now find $n$ and $\beta$ such that $Y\subseteq K(n,\beta)$. Since
$K(n,\beta)$ is closed, $Z\subseteq\overline{Y}\subseteq
K(n,\beta)$.

\noindent Now pick $M\prec H(\Theta)$ with $\Theta$ large enough
such that $|M|=\lambda\subseteq M$ and $\lambda, X,
(S_\eta:\eta<\lambda)\in M$. Then $S_\eta\cap\overline{X}\cap
S_{\kappa}^{\leq\lambda}$ is non-empty for all $\eta$. By
elementarity, $M$ sees this and so $M\cap S_\eta\cap\overline{X}\cap
S_{\kappa}^{\leq\lambda}$ is non-empty for all $\eta$. However
$M\cap\overline{X}\cap S_{\kappa}^{\leq\lambda}$ has size $\lambda$
so there are $n$ and $\beta$ such that $M\cap\overline{X}\cap
S_{\kappa}^{\leq\lambda}\subseteq K(n,\beta)$. So $S_\eta\cap
K(n,\beta)$ is non-empty for all $\eta$. \qed

\section{Strongly compact cardinals and $\CP(\kappa,\theta)$} \label{secCPLC}
\noindent We turn to the proof of theorem \ref{thSCSCH}. We will
need the following trivial consequence of the existence of a
strongly compact cardinal:
\begin{lemma}
Assume $\lambda$ is strongly compact. Then for every regular
$\kappa\geq\lambda$, there is $\U$, $\lambda$-complete uniform
ultrafilter on $\kappa$ which concentrates on $S_\kappa^{<\lambda}$.
\end{lemma}
\proof Assume $\lambda$ is strongly compact and $\kappa\geq\lambda$
is regular. By definition of $\lambda$ there is a $\lambda$-complete
ultrafilter $\W$ on $[\kappa]^{<\lambda}$ such that for all
$X\in[\kappa]^{<\lambda}$, $\{Y\in[\kappa]^{<\lambda}:X\subseteq
Y\}\in\W$. Set $\U$ to be the family of $A\subseteq\kappa$ such that
$\{X\in[\kappa]^{<\lambda}:\sup(X\cap\a)=\a\}\in\W$ for all $\a\in
A$. It is immediate to check that $\U$ is a $\lambda$-complete
ultrafilter which concentrates on $S_\kappa^{<\lambda}$. \qed

\noindent Now let $\theta<\lambda$ and $\kappa\geq\lambda$ be
regular cardinals and fix a $\theta$-covering matrix $\D$
$=(K(\a,\beta):\a\in\theta,\beta\in\kappa)$ for $\kappa$. Let
$A_\a^\gamma=\{\beta>\gamma:\gamma\in K(\a,\beta)\}$ and
$A_\a=\{\gamma\in S_\kappa^{<\lambda}:A_\a^\gamma\in\U\}$. Since
$\theta<\lambda$, by the $\lambda$-completeness of $\U$, for every
$\gamma\in S_{\kappa}^{<\lambda}$, there is a least $\a$ such that
$A_\a^\gamma\in\U$. Thus $\bigcup_{\a<\theta}
A_\a=S_{\kappa}^{<\lambda}$. So there is $\a<\theta$ such that
$A_\a\in\U$. In particular $A_\a$ is unbounded. Now let $X$ be a
subset of $A_\a$ of size $\theta$. Then $A_\a^\gamma\in\U$ for all
$\gamma\in X$. Since $|X|=\theta<\lambda$, $\bigcap_{\gamma\in X}
A_\a^\gamma\in\U$ and thus is non-empty. Pick $\beta$ in this latter
set. Then $X\subseteq K(\a,\beta)$. Since $X$ is an arbitrary subset
of $A_\a$ of size $\theta$, we conclude that $[A_\a]^\theta$ is
covered by $\D$. This concludes the proof of theorem\footnote{Notice
that for the proof of this theorem we just needed property {\it (i)}
of a covering matrix. The proof of $\CP(\kappa,\w)$ assuming either
$\PID$ or $\MRP$ will need properties {\it (i), (iii), (iv), (v)}.}
\ref{thSCSCH}.\qedd

\section{$\PID$ implies $\CP$} \label{secCPPID}

\noindent We turn to the proof of theorem \ref{thPIDSCH}. As we will
see below a model of $\PID$ retains enough properties of the
supercompact cardinals from which it is obtained in order that a
variation of the above argument can be run also in the context of
$\w$-covering matrices. We break the proof of theorem \ref{thPIDSCH}
in two parts. Assume $\kappa$ is regular and let
$\D=(K(n,\a):n\in\w,\,\a<\kappa)$ be a covering matrix on $\kappa$.
Let $\I$ be the family of $X\in[\kappa]^\w$ such that $X\cap
K(n,\a)$ is finite for all $\a< \kappa$ and for all $n<\w$.

\begin{claim} \label{clPIDCP}
$\I$ is a $P$-ideal.
\end{claim}

\proof Let $\{X_n:n\in\w\}\subseteq\I$. Let $Y=\bigcup_n X_n$. Let
$\gamma_Y$ witness {\it (v)} for $\D$ relative to $Y$. Now since for
every $n,m$, $X_n\cap K(m,\gamma_Y)$ is finite, let $X(n,m)$ be the
finite set
$$X_n\cap K(m,\gamma_Y)\setminus K(m-1,\gamma_Y)$$
\noindent and let:
$$
X=\bigcup_n\bigcup_{j\geq n} X(n,j).
$$
\noindent Notice that $X_n=\bigcup_j X(n,j)$ and $\bigcup_{j\geq
n}X(n,j)\subseteq X$, so we have that $X_n\subseteq^* X$. Moreover
$X\cap K(n,\gamma_Y)=\bigcup_{j\leq i\leq n}X(j,i)$, so it is
finite. We claim that $X\in\I$. If not there would be some $\beta$
and some $l$ such that $X\cap K(l,\beta)$ is infinite. Now $X\cap
K(l,\beta)\subseteq Y\cap K(l,\beta)\subseteq K(m,\gamma_Y)$ for
some $m$. Thus we would get that $X\cap K(m,\gamma_Y)$ is infinite
for some $m$ contradicting the very definition of $X$.\qed

\noindent Now notice that if $Z\subseteq\kappa$ is any set of
ordinals of size $\aleph_1$ and $\a=\sup(Z)$, there must be an $n$
such that $Z\cap K(n,\a)$ is uncountable. This means that
$\I\not\subseteq [Z]^\w$, since any countable subset of $Z\cap
K(n,\a)$ is not in $\I$. This forbids $\I$ to satisfy the first
alternative of the $P$-ideal dichotomy. So the second possibility
must be the case, i.e. we can split $\kappa$ in countably many sets
$A_n$ such that $\kappa=\bigcup_n A_n$ and for each $n$,
$[A_n]^\w\cap\I=\emptyset$.

\begin{claim}
For every $n$, $[A_n]^\w$ is covered by $\D$.
\end{claim}

\proof Assume that this is not the case and let $X\in[A_n]^\w$ be
such that $X\setminus X(l,\beta)$ is non-empty for all $l,\beta$.
Now let $X_0$ be a subset of $X$ such that $X_0\cap K(l,\gamma_X)$
is finite for all $l$. Then exactly as in the proof of claim
\ref{clPIDCP} we can see that $X_0\in[A_n]^\w\cap\I$. This
contradicts the definition of $A_n$.\qed

\noindent This concludes the proof of theorem \ref{thPIDSCH}.\qedd

\section{$\MRP$ implies $\SCH$} \label{secCPMRP}

\noindent We prove theorem \ref{thMRPSCH}. Thus assume $\MRP$ and
let $\D$ be a covering matrix on a regular $\kappa>\cc$ such that
$K(n,\beta)$ is a closed set of ordinals for all $n$ and $\beta$.
Assume that for all $A$ unbounded in $\kappa$, $[A]^\w$ is not
covered by $\D$. We will reach a contradiction. For each
$\delta<\kappa$ of countable cofinality, fix $C_\delta$ cofinal in
$\delta$ of order type $\w$. Let $M$ be a countable elementary
submodel of $H(\Theta)$ for some large enough regular $\Theta$. Let
$\delta_M=\sup(M\cap\kappa)$ and $\beta_M$ be the ordinal
$\gamma_{M\cap\kappa}$ provided by property {\it (v)} of $\D$
applied to $M\cap\kappa$. Set $\Sigma(M)$ to be the set of all
countable $X\subseteq M\cap\kappa$ bounded in $\delta_M$ such that
$$
\sup(X)\not\in K(|C_{\delta_M}\cap\sup(X)|,\beta_M).
$$
\noindent We will show that $\Sigma(M)$ is open and $M$-stationary.
Assume this is the case and let $\{M_\eta:\eta<\w_1\}$ be a
reflecting sequence for $\Sigma$. Let $\delta_{M_\xi}=\delta_\xi$
and $\delta=\sup_{\w_1}\delta_\xi$. Find $C\subseteq\w_1$ club such
that $\{\delta_\xi:\xi\in C\}\subseteq K(n,\delta)$ for some $n$
(which is possible since the $K(n,\delta)$ are closed subsets of
$\kappa$). Let $\a$ be a limit point of $C$. Let $M=M_\a$ and notice
that by our choice of $\beta_M$ for all $m$, there is $l$ such that
$K(m,\delta)\cap M\subseteq K(l,\beta_M)$. This means that for all
$\eta\in C\cap\a$, $\delta_\eta\in K(n,\delta)\cap M\subseteq
K(l,\beta_M)$ for some fixed $l$. Since $\a$ is a limit point of $C$
there is $\eta\in\a\cap C$ such that
$|C_{\delta_M}\cap\delta_\eta|>l$ and
$M_\eta\cap\kappa\in\Sigma(M)$. But this is impossible, since
$M_\eta\cap\kappa\in\Sigma(M)$ means that $\delta_\eta\not\in
K(|C_{\delta_M}\cap\delta_\eta|,\beta_M)$, i.e. $\delta_\eta\not\in
K(l,\beta_M)$.

\smallskip

\noindent We now show that $\Sigma_M$ is open and $M$-stationary:

\begin{claim} $\Sigma(M)$ is open.
\end{claim}
\proof Assume $X\in\Sigma(M)$, we will find $\gamma\in X$ such that
$[\{\gamma\},X]\subseteq\Sigma(M)$. To this aim notice that
$C_{\delta_M}\cap\sup(X)$ is a finite subset of $X$. Let
$n_0=|C_{\delta_M}\cap\sup(X)|$ and
$\gamma_0=\max(C_{\delta_M}\cap\sup(X))+1$. Since $X\in\Sigma(M)$,
$\sup(X)\not\in K(n_0,\beta_M)$ and so, since $K(n_0,\beta_M)$ is
closed, $\gamma_1=\max(K(n_0,\beta_M)\cap\sup(X))<\sup(X)$. Thus,
let $\gamma\in X$ be greater or equal than
$\max\{\gamma_1+1,\gamma_0\}$. If $Y\in[\{\gamma\},X]$, then
$\gamma_0\leq\sup(Y)\leq\sup(X)$, so
$|C_{\delta_M}\cap\sup(Y)|=|C_{\delta_M}\cap\sup(X)|=n_0$ and
$$
\gamma_1=\max(K(n_0,\beta_M)\cap\sup(X))<\sup(Y)\leq\sup(X)<\min(K(n_0,\beta_M)\setminus\sup(X)).
$$
\noindent Thus $Y\not\in K(|C_{\delta_M}\cap\sup(Y)|,\beta_M)$, i.e.
$Y\in\Sigma(M)$. \qed

\begin{claim}
$\Sigma(M)$ is $M$-stationary.
\end{claim}
\proof Let $f:[\kappa]^{<\w}\rightarrow\kappa$ in $M$. We need to
find $X\in\Sigma(M)$ such that $f[[X]^{<\w}]=X$. Let $N\prec
H(\kappa^+)$ be a countable submodel in $M$
 such that $f\in N$ and let $C=\{\delta<\kappa:f[[\delta]^{<\w}]=\delta\}$. Let also
 $n_0=|C_{\delta_M}\cap\sup(N\cap\kappa)|$ and $\gamma_0\in N$ be larger than
 $\max(C_{\delta_M}\cap\sup(N\cap\kappa))$. Then $(C\setminus\gamma_0)\in N$.
We assumed that no $A$ unbounded in $\kappa$ is such that $[A]^\w$
is covered by $\D$. So in particular by elementarity of $N$:
$$
N\models[(C\setminus\gamma_0)\cap S_\kappa^\w]^\w\mbox{\it is not
covered by }\D
$$
\noindent Thus there exists $X\in N$ countable subset of
$(C\setminus\gamma_0)\cap S_\kappa^\w$ such that for all $n$ and
$\beta$, $X\setminus K(n,\beta)$ is non-empty. Let $\gamma\in
X\setminus K(n_0,\beta_M)$. Now find $Z\in N$ countable and cofinal
in $\gamma$ and let $Y$ be the $f$-closure of $Z$. Then $Y\in
N\subseteq M$. Now $\gamma\in C$ so $\sup(Y)=\sup(Z)=\gamma\not\in
K(n_0,\beta_M)$. Moreover $\gamma=\sup(Y)\in
(C\setminus\gamma_0)\cap N$, so
$\gamma_0<\sup(Y)<\sup(N\cap\kappa)$, i.e.
$|C_{\delta_M}\cap\sup(Y)|=|C_{\delta_M}\cap\sup(N\cap\kappa)|=n_0$.
Thus:
$$
\sup(Y)\not\in K(|C_{\delta_M}\cap\sup(Y)|,\beta_M).
$$
\noindent I.e. $Y\in\Sigma(M)$. \qed

\noindent This concludes the proof of theorem \ref{thMRPSCH}.\qedd


\section{"Saturation" properties of models of strong forcing axioms.}

Since forcing axioms have been able to settle many of the classical
problems of set theory, we can expect that the models of a forcing
axiom are in some sense categorical. There are many ways in which
one can give a precise formulation to this concept. For example, one
can study what kind of forcing notions can preserve $\PFA$ or $\MM$,
or else if a model $V$ of a forcing axiom can have an interesting
inner model $M$ of the same forcing axiom. There are many results in
this area, some of them very recent. First of all there are results
that shows that one has to demand a certain degree of resemblance
between $V$ and $M$. For example assuming large cardinals it is
possible to use the stationary tower forcing introduced by
Woodin\footnote{\cite{larGU} gives a complete presentation of this
subject.} to produce two transitive models $M\subseteq V$ of $\PFA$
(or $\MM$ or whatever is not conflicting with large cardinal
hypothesis) with different $\w$-sequences of ordinals and an
elementary embedding between them. However $M$ and $V$ do not
compute the same way neither the ordinals of countable cofinality
nor the cardinals. On the other hand, K\"onig and Yoshinobu
{\cite[Theorem 6.1]{KY}} showed that $\PFA$ is preserved by
$\w_2$-closed forcing while it is a folklore result that $\MM$ is
preserved by $\w_2$-directed closed forcing. Notice however that all
these forcing notion do not introduce new sets of size at most
$\aleph_1$. In the other direction, in \cite{velSCH}
Veli\v{c}kovi\'c used a result of Gitik to show that if $\MM$ holds
and $M$ is an inner model such that $\w_2^M=\w_2$, then ${\cal
P}(\w)\subseteq M$ and Caicedo and Veli\v{c}kovi\'c
\cite{caivelBPFA} showed, using the mapping reflection principle
$\MRP$ introduced by Moore in \cite{mooSCH}, that if $M\subseteq V$
are models of $\BPFA$ and $\omega_2^M=\omega_2$ then ${\cal
P}(\w_1)\subseteq M$. In any case all the results so far produced
show that any two models $V\subseteq W$ of some strong forcing axiom
and with the same cardinals have the same $\w_1$-sequences of
ordinals. Thus it is tempting to conjecture that forcing axioms
produce models of set theory which are "saturated" with respect to
sets of size $\aleph_1$. One possible way to give a precise
formulation to this idea may be the following:

\begin{conjecture} \label{conrigw1} (Veli\v{c}kovi\'c) Assume $W\subseteq V$ are models
of $\MM$ with the same cardinals. Then $[Ord]^{\leq\w_1}\subseteq
W$.
\end{conjecture}
\noindent or in a weaker form:

\begin{conjecture} \label{conrigw2} Assume $W\subseteq V$ are models
of $\MM$ with the same cardinals. Then $[Ord]^{\w}\subseteq W$.
\end{conjecture}

\noindent $\CP$ can be used to show that one shouldn't expect that
the negation of conjecture \ref{conrigw2} can be obtained by means
of forcing. Since no inner model theory is known for models of
$\MM$, $\PFA$ or of a supercompact cardinal, the results I will
present brings to the conclusion that currently there are no
suitable means to try to prove the negation of the second
conjecture. Moreover I can show that if this conjecture fails first
at $\kappa$, then $\kappa$ is close to be a J\'onsson cardinal in
the smaller model and more.

\begin{theorem}\label{thPFCP}
Assume $\CP(\kappa^+,\theta)$. Let $W$ be an an inner model such
that $\kappa$ is a regular cardinal of $W$ and such that
$(\kappa^+)^W=\kappa^+$. Then $\cf(\kappa)\neq\theta$.
\end{theorem}

\noindent This shows that if $\lambda$ is strongly compact, than one
cannot change the cofinality of some regular $\kappa\geq\lambda$ to
some $\theta<\lambda$ and preserve at the same time $\kappa^+$ and
the strong-compactness of $\lambda$.

\begin{corollary} \label{corPFCP}
Assume $\PFA$ and let $W$ be an inner model with the same cardinals.
Then $W$ computes correctly all the points of countable cofinality.
\end{corollary}

\proof $\PFA$ implies $\CP$. Now apply the above theorem. \qedd

\noindent In particular this gives another proof that Prikry forcing
destroys $\PFA$, since if $g$ is a Prikry generic sequence on a
measurable $\kappa$, $V[g]$ cannot model $\CP$.

\noindent The next propositions show that neither $\kappa^+$-CC
forcing notions nor homogeneous forcing notions can be used to
obtain a generic extension which models $\PFA$ and which first adds
a new $\w$-sequence to $\kappa$. First of all we can see that the
failure of conjecture \ref{conrigw2} entails that $W$ and $V$ have a
diverging notion of stationarity:

\begin{proposition} \label{thDPFCP}
Let $W$ be an inner model of $V$ and assume that $V$ models $\CP$.
Assume that for the least $W$-cardinal $\kappa$ such that
$\kappa^\omega\setminus W$ is nonempty, $\kappa^+=(\kappa^+)^W$ and
$\kappa>\cc$. Moreover assume that $[\lambda]^\w\subseteq W$ for all
$\lambda<\kappa$. Then for every $\lambda<\kappa$ regular cardinal
of $V$, there is $S$, $W$-stationary subset of
$(S_{\kappa^+}^\lambda)^W$, which is not anymore stationary in $V$.
\end{proposition}

\noindent In particular any forcing $P$ satisfying the
$\kappa^+$-chain condition and such that $\kappa$ is the least
ordinal to which $P$ adds a new $\omega$-sequence destroys $\PFA$,
since these forcings preserves stationary subsets of $\kappa^+$.
Thus also diagonal Prikry forcing kills $\PFA$.

\noindent Recall that $(P,<_P)$ is an homogeneous notion of forcing
if for every $p,q\in P$, there are $p_0\leq_P p$ and $q_0\leq_P q$
such that $(\{r: r\leq p_0\},<_P)$ is isomorphic to $(\{r: r\leq
q_0\},<_P)$. If $P$ is an homogeneous notion of forcing it is
possible to prove by a standard induction argument that if
$\phi(a_0,\cdots,a_n)$ is a formula in the language of forcing whose
parameters are canonical names for sets in the ground model, then
either $\Vdash_P\phi(a_0,\cdots,a_n)$ or
$\Vdash_P\neg\phi(a_0,\cdots,a_n)$.

\noindent Notice that the standard cardinal preserving forcing to
add an $\w$-sequence to an uncountable cardinal like Prikry forcing
or diagonal Prikry forcing are homogeneous\footnote{It should be
expected that the stationary set preserving forcing used by Magidor
to show that $\MM$ implies there is a stationary set of bad points
in $\aleph_{\w+1}$ is not homogeneous. In any case this forcing
collapses $\aleph_{\w+1}$ to an ordinal of size $\aleph_1$ and
moreover $\CH$ holds in the extension.}.

\begin{proposition} \label{thforcingkiller}
Let $P$ be an homogeneous notion of forcing, let $G$ be a
$P$-generic filter. Assume that $\Vdash_P[\lambda]^\w\subseteq V$
for all $\lambda<\kappa$ and that $\kappa^+>\cc$ is preserved.
Assume moreover that $V[G]$ models $\CP(\kappa^+)$. Then
$\Vdash_P[\kappa]^\w\subseteq V$
\end{proposition}

\noindent $\CP(\kappa^+)$ follows both from $\PFA$ and from the
existence of a supercompact $\lambda\leq\kappa$. In particular the
above proposition shows that $\PFA$, $\MM$ and supercompact
cardinals are not preserved by all known forcing notions which add
$\w$-sequences and preserve cardinals since all such forcing notions
are homogeneous.

\noindent We can conclude that there are not much hopes to find a
forcing which first falsify conjecture \ref{conrigw2} at some
$\kappa$. By theorem \ref{thPFCP}, $\kappa$ must be of countable
cofinality. Moreover we have to seek for a $P$ which on one hand
preserves cardinals and is stationary set preserving, on the other
hand is not homogeneous and does not have the $\kappa^+$-chain
condition\footnote{In fact more than that, since the proof of
proposition \ref{thDPFCP} shows that this $P$ should have the
property of forcing the non-stationarity of at least one set in
every disjoint family in the ground model of $\kappa$-many
stationary (in the ground model) subsets of
$S_{\kappa^+}^\lambda$.}.

\smallskip

\noindent The following propositions give some ideas on how to
measure the consistency strength of the negation of the above
conjectures.

\smallskip

\noindent Say that $\kappa$ is $\E$-J\'onsson for a family $\E$ of
functions $\psi:[\kappa]^{<\w}\rightarrow\kappa$ if for every
$\psi\in\E$ there is $Y$ proper subset of $\kappa$ which is
$\psi$-closed. $\kappa$ is J\'onsson if $\E$ contains all the
relevant functions.

\noindent Any measurable cardinal is J\'onsson. Silver has
shown\footnote{See \cite{kanmagLC} for a proof.} that if
$\aleph_\w>\cc$ is J\'onsson, then it is measurable in  an inner
model. It is a major open problem whether $\aleph_\w$ can be a
J\'onsson cardinal in some model of $\ZFC$ and recently K\"onig
\cite{konMMJON} has shown that $\MM$ is consistent with $\aleph_\w$
not being J\'onsson.

\begin{proposition} \label{th4IMCP}
Let $W\subseteq V$ be models of $\ZFC$ and $\kappa>\cc$ be the least
such that $[\kappa]^\w\setminus W$ is nonempty. Assume that
$V\models\CP$, $2^\w<\kappa$ and that $(\kappa^+)^W=\kappa^+$. Then
$\kappa^+$ is $W$-J\'onsson.
\end{proposition}

\begin{proposition} \label{th5IMCP}
Let $W\subseteq V$ be models of $\ZFC$ with the same bounded subsets
of $\kappa$ and computing $\kappa^+$ the same way. Assume that
$V\models\CP$, $2^\w<\kappa$, and that $\kappa$ is the least such
that $[\kappa]^\w\not\subseteq W$. Then $\kappa$ is J\'onsson in
$W$.
\end{proposition}

\noindent We now turn to the proofs.

\noindent{\bf Proof of theorem \ref{thPFCP}:} Assume the theorem is
false. Then $\cf(\kappa)=\theta$. We need some preparation. Work in
$W$. Fix a sequence $\C=\{C_\xi:\xi<\kappa^+\}\in W$ such that for
all limit $\a$, $C_\a$ is a club subset of $\a$ of order type
$\cf(\a)$. If $\xi=\a+1$, $C_{\xi}=\{\a\}$. Consider the function
$\rho^*:[\kappa^+]^2\rightarrow\kappa$ defined from such a sequence
which was introduced in the proof of lemma \ref{thCPmax}.

\noindent Now work in $V$. Let $g:\theta\rightarrow\kappa$ be a
cofinal sequence. \noindent Set
$D(\a,\beta)=\{\xi<\a:\rho^*(\xi,\beta)\leq g(\a)\}$. The two lemmas
\ref{keylemrho1} and \ref{keylemrho2} allows to prove the analogue
of fact \ref{kayfacPFCP}:

\begin{fact} The following holds:
\begin{description}

\item[\it (i)] $\otp(D(\a,\beta))<\kappa$ for all $\a<\theta$ and $\beta<\kappa^+$,

\item[\it (ii)] for all $\gamma<\beta<\kappa^+$,there is $\a_0<\theta$ such that
$D(\a,\gamma)=D(\a,\beta)\cap\gamma$ for all $\a\geq \a_0$.

\end{description} \qed
\end{fact}

\noindent This means that
$\D=(D(\a,\beta):\a<\theta,\,\beta<\kappa^+)$ is a $\theta$-covering
matrix for $\kappa^+$ with $\beta_\D=\kappa$. By
$\CP(\kappa^+,\theta)$ there is $A$ unbounded in $\kappa^+$ such
that $[A]^\theta$ is covered by $\D$. We claim that $[A]^\kappa$ is
covered by $\D$. Suppose not, then there would be some
$\gamma\leq\beta$ such that $\otp(A\cap\gamma)>\kappa$ and
$A\cap\gamma\setminus D(\a,\beta)$ is non-empty for all $\a<\theta$.
Thus we could find an $X\in [A\cap\gamma]^\theta$ such that
$X\setminus D(\a,\beta)$ is non-empty for all $\a<\theta$. But then,
by {\it (ii)} of the above fact, $X\setminus D(\a,\beta)$ would be
non-empty for all $\a<\theta$ and $\beta<\kappa^+$. This would
contradict the fact that $[A]^\theta$ is covered by $\D$. Now let
$\gamma$ be such that $\kappa<\otp(A\cap\gamma)$. Then for some $\a$
and $\beta$ we would have $A\cap\gamma\subseteq D(\a,\beta)$ and so:
$$
\kappa<\otp(A\cap\gamma)\leq\otp(D(\a,\beta))<\kappa.
$$
\noindent This is the contradiction which proves the theorem. \qedd

\noindent {\bf Proof of proposition \ref{thDPFCP}:} By the previous
theorem $\kappa$ cannot be regular in $W$, otherwise $V$ cannot
model $\CP$. Since $\kappa$ is the least such that
$\kappa^\w\setminus W$ is non empty and is not regular in $W$, we
can conclude that $\kappa$ is in $W$ a cardinal of countable
cofinality. So there is a covering matrix $\mathcal{D}$ $\in W$ for
$(\kappa^+)^W$ with $\beta_\D=\kappa$ and such that $K(n,\beta)$ is
closed for all $n$ and $\beta$. Since $\kappa^+=(\kappa^+)^W$, $\D$
is still a covering matrix in $V$ for $\kappa^+$ with
$\beta_\D=\kappa$. Assume towards a contradiction that there is some
$\lambda<\kappa$ regular cardinal of $V$ such that every
$W$-stationary subset of $S_{\kappa^+}^{\lambda}$ is still
stationary in $V$. So fix in $W$ some $\{A_\alpha:\alpha<\kappa\}\in
W$ partition of $(S_{\kappa^+}^\lambda)^W$ in $\kappa$-many
$W$-stationary sets. Now by our assumption, this is still a family
of disjoint stationary subsets of $V$. Let for every $n,\beta$,
$D(n,\beta)$ be the set of $\alpha<\kappa$ such that $A_\alpha\cap
K(n,\beta)$ is nonempty. Now $D(n,\beta)\in W$ and
$|D(n,\beta)|\leq|K(n,\beta)|$ has size less than $\kappa$. So, by
minimality of $\kappa$, we have that $D(n,\beta)^\omega\subseteq W$,
else there would be a new $\omega$-sequence in $|D(n,\beta)|$. Apply
$\CP$ in $V$ and find $X$ unbounded in $\kappa^+$ such that
$[X]^{\w}$ is covered by $\mathcal{D}$. Exactly as in the proof of
fact \ref{faREFPFA}, we can see that
$$
[\overline{X}\cap S_{\kappa^+}^{\leq\lambda}]^\lambda
$$
\noindent is covered by $\D$. Now pick $M$ countable elementary
submodel containing all relevant information. Then
$A_{g(n)}\cap\overline{X}\cap M\cap S_{\kappa^+}^\lambda$ is
nonempty for all $n$, by elementarity of $M$. Now $M\cap
\overline{X}\cap S_{\kappa^+}^\lambda\subseteq K(n,\beta)$  for some
$n,\beta$. This means that $g\in D(n,\beta)^\omega\subseteq W$ and
we are done. \hfill$\blacksquare$\medskip

\noindent{\bf Proof of proposition \ref{thforcingkiller}}: As before
we can remark that $\kappa$ has countable cofinality in $V$,
otherwise the hypothesis  of the proposition are not compatible with
the conclusion of theorem \ref{thPFCP}. For this reason
$S_{\kappa^+}^\w$ has the same interpretation in $V$ and in $V[G]$.
Now we will use homogeneity to obtain a family $T\in V$ of
$\kappa$-many disjoint subsets of $S_{\kappa^+}^\w$ which are
stationary in $V[G]$, we will then proceed exactly as in the
previous proposition to obtain the desired conclusion. To this aim
let:
$$
\{A_{\a\beta}:\a<\kappa,\beta\in\kappa^+\}\in V
$$
\noindent be an Ulam matrix, i.e for all $\a<\kappa$ and
$\beta<\kappa^+$,
$$
A_{\a\beta}=\{\xi\in S_{\kappa^+}^\w: \phi_\xi(\beta)=\a\}
$$
\noindent where $\phi_\xi:\xi\rightarrow\kappa\in V$ is a bijection
for all $\xi\in[\kappa,\kappa^+)$.

\noindent A standard $\ZFC$ argument shows that there is a row $\a$
of the matrix such that the set of $\beta$ such that $A_{\a\beta}$
is stationary has size $\kappa^+$, moreover $A_{\a\beta}\cap
A_{\a\gamma}=\emptyset$ whenever $\beta\neq\gamma$.

\noindent Fix an $\a$ with this property in $V[G]$ and let $S$ be
the set of $\beta$ such that $V[G]$ models $A_{\a\beta}$ is
stationary. Now by the homogeneity of $P$ and the fact that
$A_{\a\beta}\in V$ we can conclude that:
$$
\Vdash_P A_{\a\beta} \mbox{ is a stationary subset of }
S_{\kappa^+}^\w
$$
\noindent for all $\beta\in S$. We can conclude that $S\in V$. Let
$T=\{S_\a:\a<\kappa\}\in V$ be a subset of $S$ of size $\kappa$. Now
proceed exactly as in the previous proposition. \qedd

\noindent{\bf Proof of propositions \ref{th4IMCP} and
\ref{th5IMCP}:} Remark that by theorem \ref{thPFCP} any $\kappa$ as
in the hypothesis of theorems \ref{th4IMCP} and \ref{th5IMCP} has
already countable cofinality in $W$. Fix $A\in V$ unbounded set
witnessing $\CP$ relative to a $\D\in W$ covering matrix for
$\kappa^+$ with $\beta_\D=\kappa$. Notice that this entails that
$[A]^\w\subseteq W$. Now let
$\psi:[\kappa^+]^{<\w}\rightarrow\kappa^+$ be any function in $W$.
Let $B$ be the $\psi$-closure of $A$. We claim that $B$ is a proper
subset of $\kappa^+$, otherwise if $X\in[\kappa^+]^\w\setminus W$ is
contained in $B$ there would be a countable subset $Y$ of $A$ such
that $X$ is contained in the $\psi$-closure $Z$ of $Y$. However
$[Z]^\w\subseteq W$, since $Z\in W$ is countable, so $X\in W$ and
this contradicts our choice of $X$. This shows that $\kappa^+$ is
$W$-J\'onsson. To show proposition \ref{th5IMCP}, set for every
$\psi:[\kappa]^{<\w}\rightarrow\kappa$, $(T_\psi,<_{T_\psi})$ to be
the well-founded partial order of bounded subsets $X$ of $\kappa$
such that $X$ is strictly contained in $\sup(X)$ and such that for
all $\a_1,\cdots,\a_n\in X$, if $\psi(\a_1,\cdots,\a_n)<\sup(X)$,
then $\psi(\a_1,\cdots,\a_n)\in X$. If $X,Y\in T_\psi$,
$X<_{T_\psi}Y$ if $Y$ end extends $X$ and $|X|<|Y|$. Now notice that
$T_\psi$ gets the same interpretation in $W$ and in $V$. If
$\overline{\psi}:[\kappa^+]^{<\w}\rightarrow\kappa^+$ is in $W$ an
arbitrary extension of $\psi$ and $A$ witness in $V$ the covering
property, $\overline{\psi}[A]\cap\kappa$ is a branch through
$T_\psi$ of type $\kappa$. Thus $V$ models that $>_{T_\psi}$ is
ill-founded. By absoluteness $W$ models that $>_{T_\psi}$ is
ill-founded. Notice that our requirement on the order $<_{T_\psi}$
entails that a witness that $>_{T_\psi}$ is ill founded is
necessarily a proper subset of $\kappa$ of type $\kappa$ closed
under $\psi$. Since this argument holds for an arbitrary $\psi\in
W$, we can argue that $W$ models that $\kappa$ is J\'onsson. This
concludes the proof of both propositions. \qedd

\noindent The results that I presented suggest to investigate the
following problem as the next step towards a solution of the above
conjectures:

\begin{problem} Is it possible to prove the analogue of corollary \ref{corPFCP} for
sequences of size $\w_1$? I.e. is it possible to prove that if $V$
models $\MM$, $W$ is an inner model and $\kappa$ is a regular
$W$-cardinal with $\kappa^+=(\kappa^+)^W$, then $\cf(\kappa)>\w_1$
in $V$?
\end{problem}

\noindent Assuming a positive answer to this problem, the next
proposition shows that conjecture \ref{conrigw1} cannot have a
negative solution by means of set-forcing.

\begin{proposition} Assume $V\subseteq W$ are models of $\MM$ with
the same ordinals of cofinality $\w$ and $\w_1$ and that $W$ is a
set forcing extension of $V$. Then $[Ord]^{\w_1}\subseteq V$.
\end{proposition}

\proof $W$ is a set-generic extension of $V$ by some $P$-generic
filter. Thus $P$ is a set and satisfies the $|P|^+$-chain condition.
Let $\kappa=|P|^+$. It is enough to show that
$[\kappa]^{\w_1}\subseteq V$. Now
$S_{\kappa}^{\leq\w_1}=(S_{\kappa}^{\leq\w_1})^W$. By the
$\kappa$-chain condition we get that every stationary subset of
$S_\kappa^\w\in V$ remains stationary in $W$. Now fix a partition
$\{S_\a:\a<\kappa\}$ in $V$ of $S_{\kappa}^{\w}$ in $\kappa$-many
disjoint stationary subset of $\kappa$. Then this is still a
partition in stationary sets of $W$. Given any
$g\in\kappa^{\leq\w_1}$ find by $\MM$ in $W$ an ordinal
$\delta<\kappa$ of cofinality $\w_1$ such that $S_{g(\xi)}$ reflects
on $\delta$ for all $\xi$ in the domain of $g$. Let $C\in V$ be a
club subset of $\delta$ of order type $\w_1$ such that $\a\in C$ iff
there is $\xi$ such that $\a\in S_{g(\xi)}$. Then $Im(g)\in V$
since:
$$
Im(g)=\{\eta:\exists\a\in S_\eta\cap C\}.
$$
\noindent We can conclude that also $g\in V$. \qedd

\noindent Thus also the following problem seems to be relevant to
our discussion:

\begin{problem} Is it possible to produce two models $W\subseteq V$ of $\MM$ with the
same cardinals, the same ordinals of cofinality $\w$ and $\w_1$ and
such that some stationary set of $(S_\kappa^\w)^W$ is not anymore
stationary in $V$, where $\kappa$ is a regular cardinal of $V$?
\end{problem}

\bibliographystyle{plain}
\bibliography{bibliothesis2MV}
\end{document}